# Beginning Mathematical Writing Assignments

Alexander Halperin, Colton Magnant, Zhuojun Magnant

October 27, 2018

**Abstract**

Writing assignments in any mathematics course always present several challenges, particularly in lower-level classes where the students are not expecting to write more than a few words at a time. Developed based on strategies from several sources, the two small writing assignments included in this paper represent a gentle introduction to the writing of mathematics and can be utilized in a variety of low-to-middle level courses in a mathematics major.

## 1 Introduction

It is often a struggle to get mathematics students to write much of anything beyond their name. In fact, many math majors specifically choose mathematics believing they can avoid writing altogether. When faced with the brutal reality that writing in mathematics is just as important and in fact far more detail-sensitive than in most other disciplines, many students seriously consider changing majors. The problem is cyclic; many algebra, precalculus, and calculus classes require little to no mathematical writing, either through problem-solving, exposition, or even self-reflection after an assignment. Students, in turn, become accustomed to this lack of writing and associate mathematics with poor (at best) explanations.

The short writing assignments contained in this work together represent our attempt at a solution to several problems at once. First, these assignments provide a gentle, mathematical baby-step into the formal writing of mathematics. Second, they were used as part of a larger research project in which the authors compared the progress of students from two different classes at two different universities, one class with an intervention of a larger writing assignment [6] and the other without. We omit the details of this larger assignment as our focus here is on the smaller writing assignments. Third and related to the previous item, these assignments, with minimal background introduction, can be assigned in virtually any course within a math program. In short, they were developed as both a measurement tool and a gentle introduction to formal mathematical writing.

These writing assignments were developed using a combination of strategies including ideas from Walvoord [11], Bahls [1], Bean [2], and Crannell et al. [4] as well as ideas stemming from our own personal experience in taking and teaching courses in mathematics over the years. The assignments were intended to foster effective writing habits and, at the same time, develop students' skills in the areas of argumentation, analysis and synthesis. In order to prove their claimed solution, the students must argue using credible evidence and supporting logic. Effective analysis and description of the situation in both questions along with their corresponding difficulties is critical to a successful complete solution. Finally, an element of synthesis is expected in the summarizing conclusion, where the students must consider a possible "natural" next step as a direction for future work.

Both of the following writing assignments were assigned in each of two different classes, one at each of two different universities. University A is a mid-sized public regional comprehensive university while University B is a large public regional comprehensive university. At University A, the assignments were used in an introductory Discrete Mathematics course. At University B, the assignments were used in a course on introduction to proofs called Mathematical Structures. Each class had about 30 students, primarily second-year undergraduates with a handful of first-year and third-year students as well. Most of these students had never written more than a sentence or two in a math class. Both classes were learning Claim-Proof form[1] of mathematical writing so an additional goal of these assignments was to reinforce this writing style.

For both of the writing assignments, students spent a day in class solving a related problem so that outside of class, they could focus almost exclusively on the writing component rather than dwelling on the mathematics. The mathematical content of the assignments was also deliberately involved — students needed to spend the entire period understanding and answering the problem — but not terribly difficult to further encourage the students to focus on the writing. Classroom discussions centered around the mathematics and a bit of outlining and formatting of the written work. Armed with at least the bulk of a solution to the mathematical content of the questions and, in most cases, a rough outline of their papers, the students proceeded home to complete the papers. Anecdotally, students suggested that they were better able to synthesize information through the writing than standard exercises, and they appreciated the experience.

---

[1] Mathematicians use "Claim-Proof" form to state and logically justify an assertion of fact.

They were then expected to complete each assignment in about a week. We required that each paper contain the four semi-standard sections found in most mathematical papers: an abstract, introduction, main results section, and conclusion. We also provided classroom discussion guidance as to what type of content was expected in each section.

In both classes (at the two different universities), the student products were scored using the same writing rubric (see the first appendix for a sample from the rubric, 3 of the 17 traits) for consistency. This rubric was developed primarily based on the well-thought-out Georgia Southern University Quality Enhancement Plan writing rubric [7] but also drawing on insights from [1]. This rubric was used since it provides a reliable and consistent measure of the different components of the written products that we wanted to measure for this project.

## 2  Assignment 1: Cat and Mouse

A cat chases a mouse in and out of a house whose floor plan is shown below. Due to the hot weather and malfunctioning air conditioner, all doors and windows are open. This provides a rousing game of tag, as both the mouse and the cat are small enough to fit through all doorways and window frames. Is is possible for the cat and mouse to run through every doorway and window frame exactly once? If so, then draw such a route. If not, then prove that such a route is not possible.

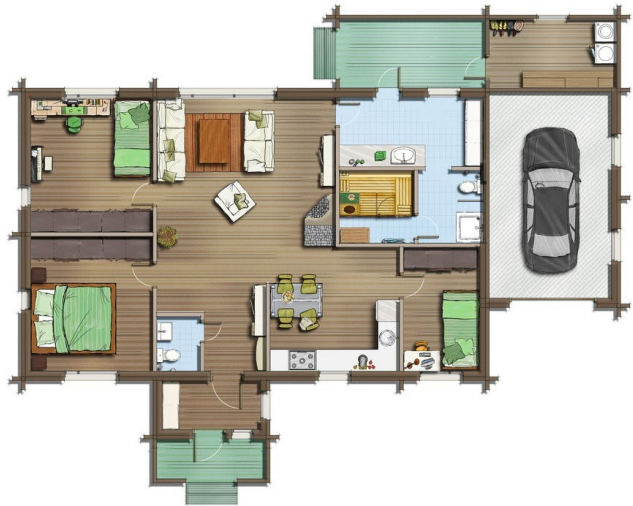

Figure 1: Image from http://zazio.xyz/maison-modernes/plan-de-maison-entunisie-100m2.html

Make sure to write up your proof in Claim-Proof form, stating the answer at the beginning with a claim and using complete sentences and paragraphs in your proof. Be sure to include any figures that may assist the reader when reading your answer. You should motivate your result and define all necessary terminology. Expect to write about a page of typed text. Your final work will be scored using the "Writing rubric" posted on the course website. Your paper should consist of the following sections:

- Abstract: briefly state the problem and the intent of your paper,

- Introduction: define relevant mathematical concepts and briefly discuss this question and how it relates to the Bridges of Königsberg problem (as we discussed in class),

- Main Result(s): state and prove your result(s),
- Conclusion: summarize your work, and make conjectures that arise from your result(s).

## 3   Assignment 2: Poker Hands

During a 5-card Poker game between three of the most famous (fictional) Poker players, tension rises when James Bond [3], Kenny Rogers [10], and Rusty Ryan [9] each go "all in," putting a combined $5 million into the pot. The situation resembles something like this:

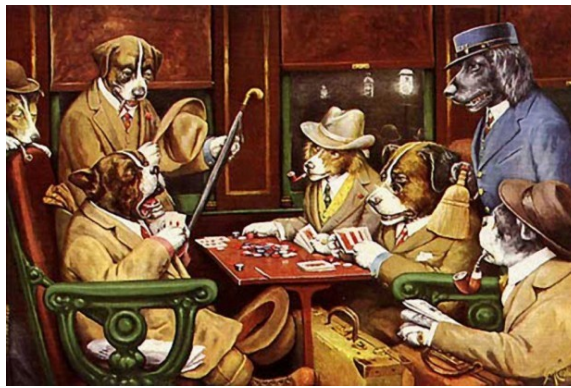

Figure 2: *His Station and Four Aces*, C. M. Coolidge, 1903.

The players reveal their hands to find that

- James Bond has a \_\_\_\_\_\_\_\_\_\_\_\_\_\_\_\_\_\_\_\_\_\_\_,
- Kenny Rogers has a \_\_\_\_\_\_\_\_\_\_\_\_\_\_\_\_\_\_\_\_\_, and
- Rusty Ryan has a \_\_\_\_\_\_\_\_\_\_\_\_\_\_\_\_\_\_\_\_\_\_.

Of course, no one wants to let go of any money. In fact, each player demands to know the exact likelihood of each hand; only then can the winner be declared. Since each player has a **different** hand, this will require three separate computations.

As the dealer, you must determine the winner. Find the general probabilities of each of the five Poker hands—that is, you must state how likely it would be to get each of the hands after drawing 5 cards from a 52-card deck (consisting of 13 values, each with 4 suits). Naturally, the hand with the lowest probability wins. It is important that you prove your answers accurately and concisely, in no more than 2 or 3 pages.

Make sure to write up your proofs in Claim-Proof form, stating the answer at the beginning with a claim and using complete sentences and paragraphs in your proof. Write a separate claim and proof for each player's hand. While you may not need any figures to assist you, you must use proper notation when referring to combinations and permutations.

Your paper should consist of the following sections:

- Abstract: briefly state the problem and the intent of your paper,

- Introduction: state the basic history and rules of Poker; also define combinations and probability,

- Main Result(s): state and prove your result(s),

- Conclusion: summarize your work, and make conjectures that arise from your result(s).

You may choose any three of the following (non crossed-out) Poker hands:

- **Royal Flush:** ~~The values 10, J, Q, K, A of the same suit.~~

- **Straight Flush:** ~~Any 5 consecutive values with the same suit.~~

- **Four of a Kind:** ~~All 4 copies of the same value and one additional card.~~

- **Full House:** Any 3 copies of one value and any 2 copies of a different value.

- **Flush:** Any 5 cards of the same suit that do not form a Royal Flush or
Straight Flush.

- **Straight:** Any 5 consecutive values that do not form a Royal Flush or a Straight Flush.

- **Three of a Kind:** Any 3 copies of one value and any 2 different values.

- **Two Pair:** Any 2 copies of one value and any 2 copies of another value and one additional value.

- **Pair:** Any 2 copies of one value and any 3 different values.

- **High Card:** All other Poker hands not previously described.

# 4 Results, Discussion and Conclusion

The playful nature [2] of both questions is intended to welcome, rather than intimidate, students when they first read the assignment in class. Immediately after reading the problem, students receive a handout with questions pertaining to the paper. For the rest of the class period, students work in groups to discover a solution. The instructor tours the classroom, answering minor questions when necessary. We detail this approach for each assignment in the subsequent sections.

Students collaborated during class but were expected to each write their own separate paper, as opposed to Latulippe and Latulippe's assignments in [14], in which 2–3-student groups turned in a single essay. This was to ensure that each student was responsible for their own writing, as we prioritized writing over problem-solving. (Although the QEP rubric weights math and writing equally, most students understood the solutions to both writing assigments by the end of the lecture, meaning that their math scores should have been high with little variance.)

## 4.1 Cat and Mouse

For the "Cat and Mouse" assignment, students first explore the famous "Seven Bridges of Königsberg problem" [8], in which a traveler tries, in a continuous route, to cross each of seven bridges in Königsberg, Prussia exactly once. After some trial and error, students discover that such a route is impossible, along with the realization that the lack of solution must be *proved*, rather than *asserted*. Listed below the problem are the steps to the proof, out of order[3], for students to rearrange. At the end of the exercise, the class recaps the argument to the instructor, in the students' own words. From there, the students have the necessary mathematical tools to answer the essay question. See the second appendix for this in-class assignment. (As a matter of fact, the "Cat and Mouse" problem has the same solution as the Seven Bridges problem: the fact that there exist more than two rooms with an odd number of doorways plus windows means the cat and mouse cannot pass through every open door and window exactly once.)

---

2 The authors' tongue-and-cheek approach was inspired by Gavin LaRose [4].
3 This idea is courtesy of Dr. Annalisa Crannell.

Given that the students can now easily solve the "Cat and Mouse" problem, the main challenge of this assignment is for students to present a mathematical solution with its proper motivation and background. By writing mostly about the framing of a problem, students are meant to see that a mathematical paper requires far more than a problem statement and solution. Additionally, they are forced to write in paragraph form, a first for most of them in a math class.

According to the rubric scores, the students particularly struggled with audience awareness by failing to provide the necessary context and definitions needed for full understanding. Broadly speaking, other areas of the rubric were in the acceptable range, especially the proof write-up. This was not surprising since the question is simple to understand and the students already have a correct proof from their class notes.

## 4.2 Poker Hands

The purpose of the poker writing assignment is to again show students that mathematics is primarily carried out in words, rather than symbols. Further, students should understand the power of the "combinatorial proof," in which quantities are counted using basic multiplication, factorials, and combinations. A problem that could take hundreds of algebraic calculations by brute force can sometimes be answered in a few short sentences in a combinatorial proof, hence making it the far more desirable option.

The writing assignment on poker hands comes at the end of a week's worth of combinatorics lessons. Students have studied the *combination* "$n$ choose $r$," or $nCr$, and learned its combinatorial definition ("$nCr$ is the number of ways to choose a subset with $r$ elements out of a set with $n$ elements") as well as derived its algebraic formula ($nCr = n!/(n-r)!r!$). Further, they have dealt with counting problems, including a poker hand example and several more involving playing cards. Most importantly, students have learned to solve these problems by viewing $nCr$ through its combinatorial definition, which emphasizes exposition and conceptual understanding over calculation. Thus, the algebraic formula is more a technical result and is only used at the end of a problem to find an exact numerical answer.

After receiving the poker hands assignment at the beginning of class, students spend the rest of the period working in groups to count the number of each type of poker hand displayed on the second page (including the crossed out hands). The professor roams the classroom, sorting out any misconceptions and correcting what are usually small errors. By the end of class, all student groups have counted all or nearly all poker hands.

The main task for their writing assignment is for students to formally write their ideas in class as a logical sequence of steps, in which they reformulate a poker hand into the values and/or suits that are chosen to form it. Of course, they still must provide background for this problem, but that is usually simpler than the previous writing assignment. Most of them find the poker hand situation more gripping than the cat-mouse chase, and many are eager to research the history of poker (and, in some cases, discuss it at great length). In fact, some students took such interest with the history of poker that it dominated the Introduction. To remedy this, future assignments will specify a five-line limit to the "history" portion of the Introduction, and the rubric will be adjusted accordingly.

This being the second assignment, the instructors were able to discuss the issues with audience awareness so the scores in this area were slightly improved over the Cat and Mouse assignment. Other areas where the students showed a bit of weakness was paragraph structure, transition between paragraphs, and overall flow. In particular, several students listed calculations for the three chosen hands with almost no discussion in between. One way to combat this misunderstanding may be to create a different rubric and provide it to students with the assignment. We discuss this at the end of the next section and give such a rubric in the Appendix.

### 4.3 General Conclusion and Future Plans

Our modest goal was to use these writing assignments simply as an introduction to mathematical writing, opening the door to the world of written mathematics. That said, given the ease in which students were able to state their solutions, we believe we should add a degree of difficulty to each assignment. For the cat and mouse assignment, this could mean asking for a more general proof that any multigraph with more than two vertices of odd degree has no Eulerian trail. For the poker assignment, we are considering requiring students to find the probability of a "high card," (the most difficult to explain), including a wild card (i.e., a card that can take on any value or suit), or perhaps a variant game with a different number of card values and/or suits.[4]

In the future, we hope to further incorporate a small reflection piece after each assignment in which the students will reflect on their process of writing, which may (should) include revision, peer-review, editing, etc. We hope to use this information to enlighten students to the effective writing processes that not only make them stronger writers but also, more immediately, result in better grades. Although explicit mention of Process writing was not included in the assignment prompts, it was discussed in class, particularly when discussing outlines, revision, and peer review. We intend to include more deliberate details about Process Writing in future iterations.

Broadly speaking, students show better performance when the audience for the assignment is clear, supporting the theory in the MAA Instructional Practices Guide [13, p. 91]. Admittedly, the specific audience was omitted from these assignment prompts since we stated in class that the intended audience was their peers, both inside and outside the course. This audience was purposely chosen to encourage peer review in the revision process. We were pleasantly surprised at the amount and effectiveness of the peer reviews and intend to be more intentional about the audience in further iterations of these assignments, in particular, including explicit statements about audience in the assignment prompts.

We also plan to develop more advanced and challenging problems to further the students' writing experience later in their mathematical careers. In both of these assignments, although students spent a day on each of the Seven Bridges of Königsberg and poker hands, the actual mathematical solutions were short, no longer than a paragraph each. Although our students appreciated the succinctness of a well-written mathematical proof, we do want them to have experience with writing a proof that, by itself, spans at least a page or, better still, requires steps or lemmas to prove.

One further change, particularly for the poker assignment, could be to require students to research a real-world combinatorial question. While similar to the approach used in [5], we would require students to find a question outside of class instead of citing previous course material. Since, unlike Pinter, our classes contain exclusively math majors and minors, we believe the added component of researching combinatorial problems outside of class would be a fair, if not challenging, additional requirement. Such an assignment would require group work for both discovering and solving the problem. To ensure that students choose an appropriately difficult problem, we would have students first hand in their problems before writing the essay.

---

[4] See *Tiny Epic Western* or *Panda´nte*.

We have also considered changing the rubric to something shorter, more readable, and with specific guidelines (e.g., "Utilize $\binom{52}{5}$ in each proof") . We believe this may reduce student forgetfulness or misunderstanding of assignment requirements, such as including Claim-Proof form for all five poker hand results or making conjectures in the conclusion. We have included an outline of such a rubric for the poker hands assignment below.

# Writing Rubric Excerpt

| Trait | Does not meet (1) | Attempted (2) | Approaches (3) | Meets (4) | Exceeds (5) |
|---|---|---|---|---|---|
| **Assignment Requirements** | The writer is off topic and/or omits most or all of the assignment requirements. | The writer addresses the appropriate topic but only superficially addresses the assignment requirements. | The writer addresses the appropriate topic and meets the assignment requirements. | The writer addresses the appropriate topic and clearly and correctly fulfills each aspect of the assignment requirements. | The writer addresses the appropriate topic and clearly, correctly, and concisely fulfills each aspect of the assignment requirements. |
| **Reasoning (proof)** | The logical connection of the argument is weak, leaving the argument or explanation unclear. A "proof by example" falls here. | The reasoning offers apparent support for the argument, but the argument or explanation is weak. | Collectively, the logic offers adequate support for the argument, but the argument or explanation remains unclear or incomplete. | Collectively, the logic supports and advances the argument or explanation of the proof. | Collectively, the logical steps offer compelling support which clearly advances the argument or explanation of the proof. |
| **Quality of Details** | Details are superficial or do not develop the proof. | Details are loosely related to the proof. Many do not provide supporting | Details are related to the proof but inconsistently provide supporting statements, credible | Details provide supporting statements, credible evidence, or the examples necessary to explain | Compelling details provide supporting statements, credible evidence, or the examples necessary to explain |

|   |   | statements, credible evidence, or the examples necessary to explain or persuade adequately. | evidence, or the examples necessary to explain or persuade adequately. | or persuade adequately. | or persuade effectively. |

| **Math 210: Discrete Mathematics** | | | **Fall 2017** |
| --- | --- | --- | --- |
| | Lecture 1: The Königsburg Bridge Problem | | Date: |
| Instructor: | | | 8/28/17 |

**Problem 1.** *Consider the following layout in Königsburg, Prussia:*

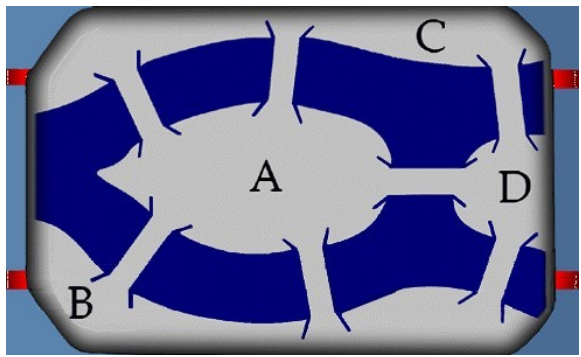

*Is there a route through the city that crosses each bridge* exactly *once?*

There are two ways we can try to figure out the answer:

1. If we think such a route exists, then ___ .

2. If we think such a route does not exist, then ___

Which method is easier? ___   Does that mean that method is correct?

Using the maps on the next page, try the easier method for a few minutes, and see what you get.

Try several different routes until you get a solution, or until you think a solution doesn't exist.

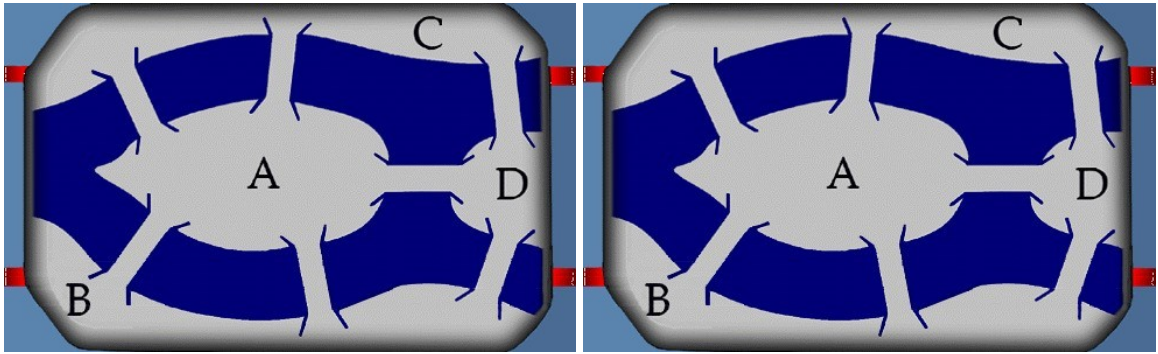

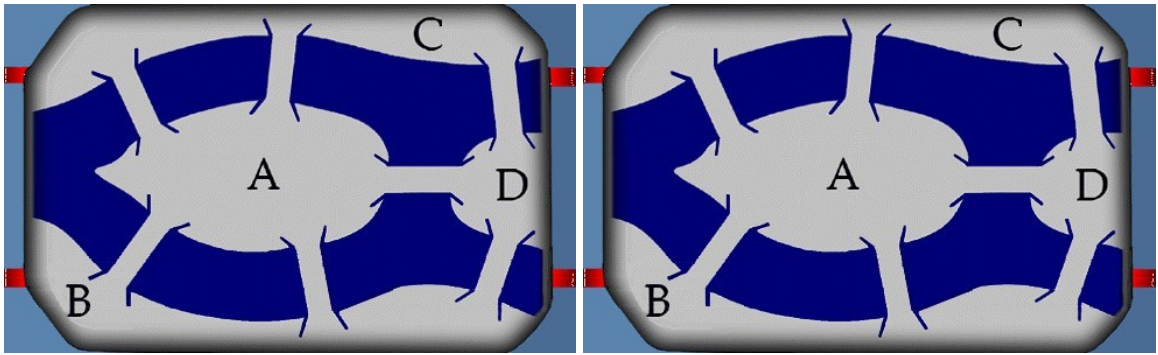

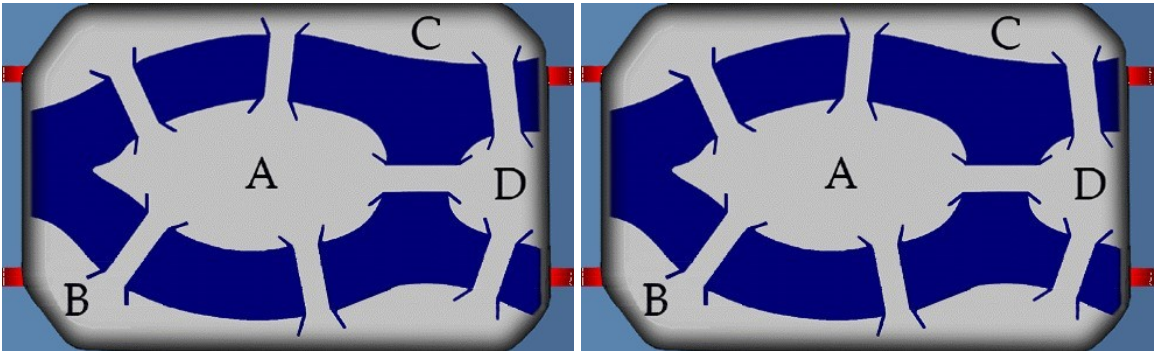

**Answer to Problem 1:** ______. (Euler, 1736)

To see this, we can view the map of K¨onigsburg as a *graph*. Each land mass is represented by a *vertex*, and each bridge is represented by an *edge*. A route through the graph that uses an edge at most one time is called a *trail*.

The steps of the *claim* and *proof* are posted below, but in the wrong order. Unscramble the steps to obtain a complete and correct claim and proof of the problem.

1. **Proof:**

2. Call this graph $G$. Hence, we see that $G$ consists of __________ vertices and __________ edges.

3. **Claim:**

4. There is no route through K¨onigsburg that traverses every bridge exactly once.

5. As a result, there cannot exist a trail in $G$ that contains every edge of $G$.

6. It now suffices to prove the claim, "There is no *trail* in $G$ that contains every *edge* of $G$."

7. However, each vertex in $G$ touches an __________ number of edges.

8. Hence, if $T$ contains every edge in $G$, then $T$ must have at least two vertices that touch an even number of edges in $T$.

9. ☐

10. Except for possibly the beginning and ending vertices, every vertex in a trail *T* touches an (circle one) *even/odd* number of edges in *T*.

11. We represent the map of Königsburg with a *graph* in the following way: draw a *vertex* for each land mass and an *edge* for each bridge.

12. This is because each middle vertex in *T* is entered by one edge and then exited by another.

Using the above reasoning, we can generalize the answer to Problem 1 as follows:

**Proposition** 2.

\_\_\_\_\_\_\_\_\_\_\_\_\_\_\_\_\_\_\_\_\_\_\_\_\_\_\_\_\_\_\_\_\_\_\_\_\_\_\_\_\_\_\_\_\_\_\_\_\_\_\_\_\_\_

\_\_\_\_\_\_\_\_\_\_\_\_\_\_\_\_\_\_\_\_\_\_\_\_\_\_\_\_\_\_\_\_\_\_\_\_\_\_\_\_\_\_\_\_\_\_\_\_\_\_\_\_\_\_

Using Proposition 2, which of the following graphs might contain a trail with every edge? For those that cannot, can you succinctly explain why? For those that might, can you find such a trail?

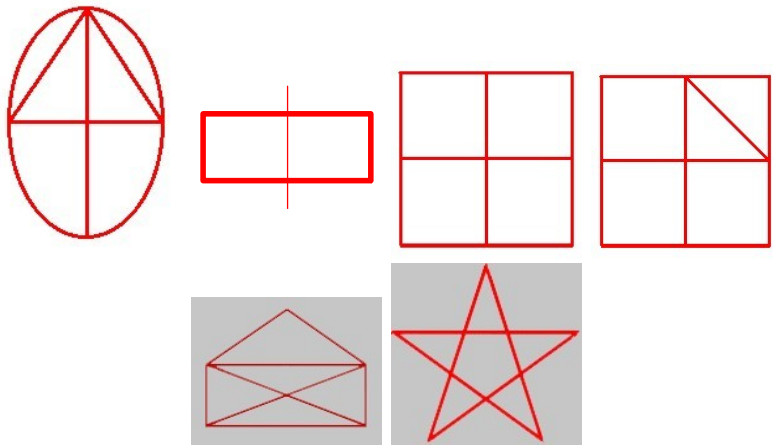

Were there any graphs where you expected to find a trail with every edge but didn't? . What does this suggest?

**Conjecture 3.**

Combining Proposition 2 and Conjecture 3, can we make an even stronger conjecture?

**Conjecture 4.**

We'll revisit this topic later in the semester...

**Homework (due Monday):** Writing Assignment #1

# Poker Paper Grading Rubric

Note that all four criteria in the Main Results section apply to each of your five main results, hence the "×5" in the Score section.

| Section | Standard | Score |
|---|---|---|
| **Abstract** | Restate the problem | /5 |
| | State the paper objective | /4 |
| | State problem-solving methods used | /1 |
| **Introduction** | Provide a brief history of poker, at most 5 lines | /10 |
| | Describe the rules of poker | /10 |
| | Restate player hands | /10 |
| **Main Results** | Use Claim-Proof form | /2 × 5 |
| | Accurately find probability | /3 × 5 |
| | Write clearly and correctly | /4 × 5 |
| | Utilize $\binom{52}{5}$ | /1 × 5 |
| **Conclusion** | Summarize results | /4 |
| | State a new question | /4 |
| | State another new question | /2 |

**Final Score:        /100**